\documentstyle[12pt]{article}
\newtheorem{proposition}{Proposition}
\newtheorem{definition}{Definition}

\newtheorem{remark}{Remark}
\newcommand{\ren}{\rho^{{\rm End}}_q}

\newcommand{\ep}{\epsilon}
\newcommand{\vko}{V_{k\omega}}
\newcommand{\vkoq}{V_{k\omega}^q}

\newcommand{\sss}{SL(n)/S(L(n-1)\times L(1))}
\newcommand{\brot}{\{\,\,,\,\,\}_{1,2}}
\newcommand{\brt}{\{\,\,,\,\,\}_2}
\newcommand{\bro}{\{\,\,,\,\,\}_1}
\newcommand{\tS}{\widetilde{S}}
\newcommand{\ahq}{A_{h,q}}

\newcommand{\uqsn}{U_q(sl(n))}

\newcommand{\vv}{V^{\otimes 2}}
\newcommand{\ww}{W^{\otimes 2}}
\newcommand{\uqs}{U_q(sl(2))}
\newcommand{\uqsl}{U_q(sl(n))}
\newcommand{\uog}{U(\overline{g})}

\newcommand{\uq}{U_q(g)}

\newcommand{\iqc}{I^q_{c_0,c_1}}
\newcommand{\om}{\omega}
\newcommand{\Co}{{\cal O}_{\om}}
\newcommand{\Cx}{{\cal O}_{x}}
\newcommand{\ogg}{\overline{g}}

\newcommand{\rd}{\rho^{\otimes 2}}
\newcommand{\vbq}{V_{\beta}^q}

\newcommand{\vb}{V_{\beta}}
\newcommand{\vqb}{V_{\beta}^q}

\newcommand{\aq}{A_{0,q}}
\newcommand{\ah}{A_{h,0}}
\newcommand{\ot}{\otimes}
\newcommand{\De}{\Delta}
\newcommand{\qq}{q+q^{-1}}

\title{Two types of Poisson pencils and related quantum objects}
\author{D.~Gurevich\\
ISTV, Universit\'e de Valenciennes,
59304 Valenciennes, France\\
J.Donin\\
Departement of Mathematics, Bar-Ilan University, 52900\\ 
Ramat-Gan,Israel\\
V.Rubtsov\\
ITEP, Bol.Tcheremushkinskaya 25, 117259 Moscow, Russia}

\date{}

\begin{document}
\maketitle
\begin{flushright}
{\it Dedicated to Alain Guichardet with regards and friendship.}
\end{flushright}

\bigskip

\begin{abstract}
Two types of Poisson pencils connected to classical R-matrices and their
quantum counterparts are considered. A representation theory of the
quantum algebras related to some symmetric orbits in $sl(n)^*$ is
constructed. A twisted version of quantum mechanics is discussed.

\end{abstract}

\section{Introduction}

There are (at least) two reasons for Poisson pencils (P.p.) to be currently
of great interest. First, they play a very important role in the theory of
integrable systems (in the so-called Magri--Lenart scheme). Second, they
appear as infinitesimal (quasi-classical) objects in the construction of
certain quantum homogeneous spaces. Roughly speaking, we can say that the
P.p. arising in the framework of the latter construction are of the first
type and those connected to integrable systems belong to the second type
(while the integrable systems themselves are disregarded).

The main characteristic that joins together these two classes of P.p. is that
a classical R-matrix participates in their construction. Let us recall that
by a classical R-matrix on a simple Lie algebra $g$ one means a
skew-symmetric $(R\in \wedge^2(g))$ solution of the classical
Yang--Baxter equation
$$[R^{12},R^{13}]+[R^{12},R^{23}]+[R^{13},R^{23}]=a\varphi,\, a\in k,$$
where $\varphi$ is a unique (up to a factor) ad-invariant element belonging
to $\wedge^3(g))$.
In what follows we deal with the ``canonical'' classical R-matrix
\begin{equation}
R=\frac12\sum X_{\alpha}\wedge X_{-\alpha}
\end{equation}
where $\{H_{\alpha},\, X_{\alpha},\, X_{-\alpha}\}$ is a
Cartan--Weyl--Chevalley base of a simple Lie algebra $g$ over the field
$k={\bf C}$. We assume that a triangular decomposition of the Lie algebra
$g$ is fixed. The field $k={\bf R}$ is also admitted, but in this case
we consider the normal real forms of the Lie algebras corresponding 
to this triangular decomposition.

Thus, the R-matrix (1) enters  the constructions of both type of P.p. under
consideration. However, the constructions and properties of these types are
completely different. Moreover, the methods of quantizing these two types of
P.p. are completely different, as well as the properties of the resulting
associative algebras. We will use the term {\em associative pencils} (a.p.)
of the first (second) type for the families of these quantum algebras arising
from the first (second) type P.p.

The ground object of the second type P.p. is a quadratic
Poisson bracket (P.b.) of Sklyanin type. By this type
we mean either the famous Sklyanin bracket\footnote{Let us specify that
this bracket is given by
$$\{\,\,,\,\,\}_S=\{\,\,,\,\,\}_l-\{\,\,,\,\,\}_r,\;\{f,g\}_{\ep}=
\mu\langle\rho_{\ep}^{\ot 2}(R),df\ot dg\rangle,\,\ep=l,r,$$
where $\rho_{l}\,(\rho_{r}):g\to {\rm Vect}({\rm Mat}(n))$ is a
representation of the Lie algebra $g$ into the space of right (left)
invariant vector fields on ${\rm Mat}(n)$, $R$ is the classical R-matrix (1)
corresponding to $sl(n)$ and $\langle\,\,,\,\,\rangle$ is the pairing between
vector fields and differentials extended to the their tensor powers.
Hereafter $\mu$ is the multiplication in the algebra under consideration.
(Let us note that for other simple Lie groups $G$ analogous P.b. being
extended to ${\rm Mat}(n)$ are no longer Poisson.)} or its various analogues:
the elliptic Sklyanin algebras in the sense of \cite{S}, the so-called
second Gelfand--Dikii structures, etc. (Although we restrict ourselves to
finite-dimensional Poisson varieties, we would like to note it is not
reasonable to consider these structures as the infinite-dimensional analogues
of P.p. on symmetric orbits in $g^*$, as it was suggested in \cite{KRR}.
The principal aim of the present paper is to make the difference between the
two classes of P.p. related to classical R-matrices more transparent.)

All these P.b. (further denoted by $\{\,\,,\,\,\}_2$) are quadratic and
their linearization gives rise to another (linear) P.b. (denoted by $\{\,\,,
\,\,\}_1$) compatible with the initial one. As a result we have a P.p.
\begin{equation}
\{\,\,,\,\,\}_{a,b}=a\{\,\,,\,\,\}_1+b\{\,\,,\,\,\}_2
\end{equation}
generated by these two brackets.

Now let us describe the corresponding quantum objects. We construct them in
two steps. First, we quantize the Sklyanin bracket and get the well-known
quadratic ``RTT=TTR'' algebra. Then by linearizing the determining quadratic
relations (or more precisely, by applying a shift operator to the algebra, cf.
below) we get a (second type) a.p. which is the quantum counterpart of the
whole P.p. These second type P.p. and their quantum counterparts are
considered in Section 2.

Let us note that the shift operator mentioned above does not give rise to any
meaningful deformation and it is reduced to a mere change of base.

We are interested in P.p. and their quantum counterparts connected to
symmetric orbits in $g^*$, where $g$ is a simple Lie algebra. It
should be pointed out that the latter P.p. are not any reduction of the
second type P.p. Let us say a little bit more about the origin of these P.p.
on symmetric orbits.

It is well known that any orbit $\Cx$ of any semisimple element in $g^*$
can be equipped with the reduced Sklyanin (sometimes called
Sklyanin--Drinfeld) bracket. One usually considers the reduction procedure
for the compact forms of simple Lie algebras, but it is also valid for normal
ones (the R-matrix for compact forms differs from (1) by the factor
$\sqrt{-1}$).

It is worth noting that the symmetric orbits in $g^*$ admit a complementary
nice property: the reduced Sklyanin bracket is compatible with the
Kirillov--Kostant--Souriau (KKS) bracket. This fact was shown in
\cite{KRR}, \cite{DG2} (in \cite{KRR} it was also shown that only symmetric
orbits possess this property).

Moreover, as shown in \cite{DG2}, both components $\{\,\,,\,\,\}_{\ep},\,\ep=
l,\,r$ of the Sklyanin bracket (cf. footnote 1) become Poisson after being
reduced to $\Cx$ and one of them (say, $\{\,\,,\,\,\}_{r}$ coincides up to
factor with the KKS one if we identify $\Cx$ with right coset $G/H$).

Since the brackets $\{\,\,,\,\,\}_{\ep},\,\ep=l,r$ are always compatible, we
have a (first type) P.p. (2) on a symmetric orbit $\Cx\in {g}^*$ with $\{\,
\,,\,\,\}_{1}=\{\,\,,\,\,\}_l^{{\rm red}}$ and $\{\,\,,\,\,\}_{2}=\{\,\,,\,\,
\}_{r}^{{\rm red}}$ (``red'' means reduced on $\Cx$)\footnote{Let us note
that similar (first type) P.p. exists on certain nilpotent orbits in $g
^*$. These pencils are generated by the KKS bracket and the so-called
R-matrix bracket defined by
$$\{\,\,,\,\,\}_{R}=\mu\langle\rho^{\ot 2}(R),df\ot dg\rangle,\;\rho={\rm ad}^*.$$
In fact the R-matrix bracket is just that $\{\,\,,\,\,\}_{l}^{{\rm red}}$ for
any such orbit but only for the symmetric ones the brackets KKS and $\{\,\,,
\,\,\}_{r}^{{\rm red}}$ are proportional to each other. The reader is
referred to \cite{GP} where all orbits in $g^*$ possessing the above P.p.
are classified.}.

The difference between these two types of P.p. also manifests itself on the
quantum level. One usually represents quantum analogues of the Sklyanin
bracket reduced to a semisimple orbit in terms of ``quantum reduction'' by
means of pairs of quantum groups or their (restricted) dual objects. However,
in this way it is not possible to represent the whole a.p. arising from the
first type P.p. under consideration. We discuss here another, more explicit,
way to describe the corresponding quantum a.p. Namely, we represent them as
certain quotient algebras.

In virtue of the results of the paper \cite{DS1}, a formal quantization of
the P.p. under consideration exists on any symmetric orbit. However, it is not
so easy to describe the quantum objects explicitly or, in other words, to
find systems of equations defining them. We consider here two ways to look
for such a system in the case when the orbit $\Cx$ is a rank 1 symmetric
space in $sl(n)^*$, i.e., that of $\sss$ type.

The first way was suggested in \cite{DG1} (but there the system of
equations was given in an inconsistent form). It uses the fact that such an
orbit can be described by a system of quadratic equations (cf. Section 3).
The second way consists in an attempt to construct a representation theory of
the first type quantum algebras (cf. Section 4) and to compute the desired
relations in the modules over these algebras. This approach is valid in
principle for any symmetric orbit in $g^*$ for any simple Lie algebra
$g$. This enables us to treat this type of quantum algebras from the
point of view of twisted quantum mechanics (Section 5). This means that
quantum algebras and their representations are objects of a twisted (braided)
category. (We consider the term ``twisted'' as more general, keeping the
term ``braided'' for nonsymmetric categories).

Completing the Introduction, we would like to make two remarks. First, let us
note that there is a number of papers devoted to q-analogues of special
functions of mathematical physics. In fact they deal with quantum analogues
of double cosets. Meanwhile, the problem of finding an explicit description
of quantum symmetric spaces (in particular, orbits in $g^*$) is
disregarded (in some sense the latter problem is more complicated). If the
usual q-special functions arise from a one-parameter deformation of the
classical objects, our approach enables us to consider their analogues
arising from the two-parameter deformation. This approach will be developed
in a joint paper of one of the author (D.G.) and L.Vainerman.

Second, we would like to emphasize that we consider the present paper as an
intermediate review of the papers \cite{DG1}, \cite{DGR}, \cite{G2} where the
present topic (i.e., the problem of the explicit description of a.p. arising
from quantization of certain P.p. associated with the ``canonical'' classical
R-matrices) was discussed. Meanwhile, this topic is still in progress. 

In what follows $\uq$-Mod stands for the category of all
finite-di\-men\-si\-o\-nal $\uq$-modules and their inductive limits. The
parameter $q$ is assumed to be generic (or formal, when we speak about
flatness of deformation, in fact we do not distinguish these two meanings).

We dedicate this paper to our friend, professor Alain Guichardet. We
acknowledge his warm hospitality at Ecole Polytechnique for a long time and
benefited a lot from our numerous discussions and conversations.

This paper was finished during the visit of V.R. to the University Lille-1.
He is thankful to his colleagues from the group of Mathematical Physics for
their hospitality. The work of V.R. was supported partially by the grants
RFFI-O1-01011, INTAS-93-2494 and INTAS-1010-CT93-0023.

\section{Second type Poisson and associative pencils}

Let us fix a simple algebra $g=sl(n),\,n\geq 2$ and consider the
corresponding quadratic Skyanin bracket extended on the space ${\rm Fun}({\rm
Mat}(n))$ which is the algebra of polynomials on the matrix elements $a_i^j$.
More precisely, this bracket is determined by the following multiplication
table
$$
\{a_{k}^{i},a_{k}^{j}\}_{2}=a_{k}^{i}\,a_{k}^{j},\,\,
\{a_{i}^{k},a_{j}^{k}\}_{2}=a_{i}^{k}\,a_{j}^{k},\,\,i<j;
$$
$$
\{a_{i}^{l},a_{k}^{j}\}_{2}=0,\,\,
\{a_{i}^{j},a_{k}^{l}\}_{2}=2\,a_{i}^{l}\,a_{k}^{j},
\,\,i<k,\,j<l.
$$

Let us linearize the bracket $\brt$. To do this we introduce a shift operator
${\rm Sh}_{h}$ which sends $a_{i}^{j}$ to $a_{i}^{j}+h \delta_{i}^{j}$ (we
extend this operator onto the whole algebra ${\rm Fun}({\rm Mat}(n))$ by
multiplicativity). Then applying this operator to r.h.s. of the above
formulas and taking the linear terms in $h$, we get a linear bracket with by
the following multiplication table
$$
\{a_{k}^{i},a_{k}^{j}\}_{1}=\delta_{k}^{i}\,a_{k}^{j}+a_{k}^{i}\,
\delta_{k}^{j},\,\,
\{a_{i}^{k},a_{j}^{k}\}_{1}=\delta_{i}^{k}\,a_{j}^{k}+a_{i}^{k}\,
\delta_{j}^{k},\,\,i<j;
$$
$$
\{a_{i}^{l},a_{k}^{j}\}_{1}=0,\,\,
\{a_{i}^{j},a_{k}^{l}\}_{1}=2\,(\delta_{i}^{l}\,a_{k}^{j}+a_{i}^{l}\,
\delta_{k}^{j}),
\,\,i<k,\,j<l.
$$
\begin{proposition} The brackets $\brot$ are compatible.
\end{proposition}

{\it Proof.} This follows immediately from the following fact: the r.h.s. of
the formulas for the bracket $\brt$ does not contain any summand of the form
$a_i ^i a_j^j$. This implies that the images of r.h.s. elements under the
operator ${\rm Sh}_{h}$ do not contain any term quadratic in $h$. The details
are left to the reader. \ $\Box$

\begin{remark}{\em It is often more convenient to take, instead of $\{a_i^j
\}$, the base consisting of $a_i^j,\,i\not=j,\,\,a_1^1-a_i^i,\,i>1$ and $a_0=
\sum a_i^i$. Then the mentioned property can be reformulated as follows: the
r.h.s. of the above formulas do not contain the term $(a_0)^2$ (with respect
to these new generators the operator ${\rm Sh}_h$ acts nontrivially only on
$a_0$).}\end{remark}

It is easy to see that the bracket $\bro$ can be represented in the form
$$
\{a_{i}^{j},a_{k}^{l}\}_{1}=\{{\bf R}(a_{i}^{j}), a_{k}^{l}\}_{gl}+
\{a_{i}^{j}, {\bf R}(a_{k}^{l})\}_{gl}.
$$
Here $\{\,\,,\,\,\}_{gl}$ is the linear bracket corresponding to the Lie
algebra $gl(n)$ (namely, $\{a_{i}^{j},a_{k}^{l}\}_{gl}=a_{i}^{l}\delta_{k}^
{j}-a_{k}^{j}\delta_{i}^{l}$) and ${\bf R}:W\to W$ is an operator defined in
the space $W={\rm Span}(a_i^j)$ as follows ${\bf R}(a_{i}^{j})={\rm sign}(j-
i)\,a_{i}^{j}$ (we assume that ${\rm sign}(0)=0$).

Thus, the space $W$ is equipped with two Lie algebra structures. The first
one is $gl(n)$ and the second one corresponds to the Poisson bracket $\{\,\,,
\,\,\}_{1}$ and therefore we have a double Lie algebra structure on the space
$W$ according to the terminology of \cite{S-T} (our construction coincides
with \cite{S-T} up to a factor).

Let us consider now the quantum analogue of the above P.p.

Let $U_q(g)$ be the Drinfeld-Jimbo quantum group corresponding to $g$
and ${\cal R}$ be the universal quantum R-matrix. Consider the linear space $V$
of vector fundamental representation $\rho : g\to {\rm End}(V)$ of the
initial Lie algebra $g$. Denote by $\rho_q$ the corresponding
representation of the quantum group $\uq$ to ${\rm End}(V)$ (we assume that
$\rho_q\to \rho$ when $q\to 1$). Then the operator $S:\vv\to\vv$ defined as
follows $S=\sigma\rho^{\ot 2}({\cal R})$ where $\sigma$ is the flip $(\sigma
(x\ot y)=y\ot x)$ satisfies the quantum Yang--Baxter equation (QYBE)
\begin{equation}
S^{12}S^{23}S^{12}=S^{23}S^{12}S^{23},\,\, S^{12}=S\ot {\rm id},\, S^{23}={\rm id}\ot S
\end{equation}
and possesses two eigenvalues (we recall that $g=sl(n)$).

Let us identify now the space $W$ with $V\ot V^*$ and equip it with the
operator $S_W=S\ot(S^*)^{-1}:\ww\to\ww$, where $S^*$ is defined with respect
to the pairing
$$
\langle x\otimes y, a\otimes b\rangle =\langle x,a\rangle\langle y,b\rangle,
\,\,x,y\in V,\,\,a,b\in V^{*}.
$$

Let us fix a base $\{a_i\}\in V$. Let $\{a^i\}\in V^*$ be the dual base and
set $a_i^j=a_i\otimes a^j$. Then we have the following explicit form for the
operator $S_W$:
$$
S_{W}(a_{i}^{k}\ot a_{j}^{l})=S^{mn}_{ij}{S^{-1}}^{kl}_{pq}
(a_{m}^{p}\ot a_{n}^{q}) \,\,{\rm where}\,\,
S(a_{i}\ot a_{j})=S_{ij}^{kl}a_{k}\ot a_{l}.
$$

It is easy to see that the operator $S_W$ also satisfies the QYBE and
possesses 1 as an eigenvalue. Then it is natural to introduce deformed
analogues of symmetric and skew-symmetric subspaces of the space $\ww$ as
follows
$$I_{-}^{q}={\rm Im}(S_{W}-{\rm id}),\,\,I_{+}^{q}={\rm Ker}(S_{W}-{\rm id}).$$
Let us describe explicitly the operator $S$ and the spaces $I_{\pm}$:
$$
S(a_i\otimes a_j)=(q-1)\delta_{i,j}a_i\otimes a_j+a_j\otimes a_i+
\sum_{i<j}(q-q^{-1})a_i\otimes a_j,
$$
$$
I_-^q={\rm Span}(a_{k}^{i}a_{k}^{j}- qa_{k}^{j}a_{k}^{i},\,
a_{i}^{k}a_{j}^{k}-qa_{j}^{k}a_{i}^{k},\, i<j;
$$
$$
a_{i}^{l}a_{k}^{j}-a_{k}^{j}a_{i}^{l},\,a_{i}^{j}a_{k}^{l}-a_{k}^{l}a_{i}^{j}
-(q-q^{-1})a_{k}^{j}a_{i}^{l},\, i<k,j<l)
$$
and
$$
I_+^q={\rm Span}((a_i^k)^2,\,q a_{k}^{i}a_{k}^{j}+a_{k}^{j}a_{k}^{i},\,
qa_{i}^{k}a_{j}^{k}+a_{j}^{k}a_{i}^{k},\, i<j;
$$
$$
a_{i}^{j}a_{k}^{l}+a_{k}^{l}a_{i}^{j},\,
a_{i}^{l}a_{k}^{j}+a_{k}^{j}a_{i}^{l}+(q-q^{-1})a_{i}^{j}a_{k}^{l},\,
i<k,j<l).
$$

Let us introduce also the quotient algebra $A_{0,q}=T(W)/\{I_{-}^{q}\}$.
In what follows $T(W)$ denotes the free tensor algebra of the space $W$
and $\{I\}$ denotes the two-sided ideal generated by the set $I\subset T(W)$
(here $I=I_-^q$ is a subspace of the space $\ww$).

It is well known (cf. \cite{RTF}) that the quantum analogue ${\rm Fun}_q(SL
(n))$ of the space ${\rm Fun}(SL(n))$ is the quotient algebra of $A_
{0,q}$ over the ideal generated by the element $\det_q-1$ ($\det_q$ is the
so-called ``quantum determinant''). Moreover, the latter quotient algebra can
be endowed with a Hopf structure. However, we are interested rather in the
algebra $\aq$ itself. This algebra is quadratic and is a flat deformation
of the classical counterpart, namely, of the symmetric algebra ${\rm Sym}(W)$
of the space $W$ (cf. \cite{DS2} for details). Let us note that the
skew-symmetric algebra of $W$ has also an evident quantum analogue
$T(W)/\{I_+^q\}$.

Now we want to introduce two parameter deformation of the algebra ${\rm Sym}
(W)$ quantizing the whole P.p. under consideration. This deformation can be
realized by means of the same shift operator as above. Applying this operator
to $I^q_-$ we obtain the following space
$$
J_{h,q}={\rm Span}(a_{k}^{i}a_{k}^{j}- qa_{k}^{j}a_{k}^{i}-
h(\delta_{k}^{i}\,a_{k}^{j}+a_{k}^{i}\,\delta_{k}^{j}),$$
$$
a_{i}^{k}a_{j}^{k}-qa_{j}^{k}a_{i}^{k}-
h(\delta_{i}^{k}\,a_{j}^{k}+a_{i}^{k}\,\delta_{j}^{k}),\,\,i<j;\;
a_{i}^{l}a_{k}^{j}-a_{k}^{j}a_{i}^{l},$$
$$
a_{i}^{j}a_{k}^{l}-a_{k}^{l}a_{i}^{j}-(q-q^{-1})a_{k}^{j}a_{i}^{l}-
hm\,(\delta_{i}^{l}\,a_{k}^{j}+a_{i}^{l}\,\delta_{k}^{j}),\, i<k,j<l),
$$
where $m=1+q^{-1}$ (we have realized here a substitution $h(q-1)\to h$).

Let us introduce the algebra $A_{h,q}=T(W)/\{J_{h,q}\}$. It is evident that
this algebra is a flat deformation of the initial commutative algebra $A_{0,
1}$ since the passage from the algebra $\aq$ to that $\ahq$ is trivial from
the deformation point of view and reduces to a change of a base. Thus, we
have constructed the a.p. $\ahq$, which is a quantum analogue of the P.p.
under consideration (it is easy to see that quasi-classical term of this two
parameter deformation is just the above second type P.p.).

A similar method can be applied to quantize the P.p. generated by the
elliptic Sklyanin P.b \cite{S}. This bracket (denoted also by
$\{\,\,,\,\,\}_2$) is determined by the following multiplication table
$$
\{S_{1}\,S_{0}\} = 2J_{23}S_{2}S_{3},\;\, \{S_{1}\,S_{2}\} = -2S_{0}S_{1}
$$
and their cyclic permutations with respect to the indices $(1,2,3)$ with some
elliptic functions $J_{ij}$ (cf. \cite{S}). The linearization of this bracket
defined by the shift operator $S_0\to S_0+h,\, S_i\to S_i, i\not=0$ ($S_0$
plays here the role of the element $a_0$ from Remark 1) gives rise to the
linear one corresponding to the Lie algebra $g=so(3)\oplus k$. The
brackets are compatible (by the same reason as above).

To quantize the P.p. generated by them, it suffices to quantize the bracket
$\brt$ and apply the shift operator. It is well known (cf. \cite{S}) that a
quantum analogue of the bracket $\brt$ is the quotient algebra $T(W)/\{I\}$
where $W={\rm Span}(S_0,....,S_3)$ and $\{I\}$ is two-sided ideal in $T(W)$
generated by the elements
$$
S_1 S_0-S_0 S_1+i J_{23}(S_2 S_3 +S_3 S_2),\,\,
S_1 S_2 -S_2 S_1 -i(S_0 S_3 + S_3 S_0)
$$
and their cyclic permutations with some elliptic functions $J_{ij}$.

We leave to the reader the explicit description of the resulting two
parameter quantum a.p.

The algebra $T(W)/\{I\}$ originally defined by E. Sklyanin is an elliptic
analogue of the symmetric algebra of the space $W$. Unfortunately, we do
not know any natural elliptic analogue of the skew-symmetric one (it seems
very plausible that such analogue does not exist).

\begin{remark}{\em There are other quantum analogues of symmetric and
skew-symmetric algebras of the space $W$ defined by the so-called reflection
equations (RE)
$$
Su_1Su_1=u_1Su_1S,\; u_1=u\ot 1, \;u=(u_i^j)
$$
where $S$ is a solution of (3). Conjecturally this ``RE algebra'' is a flat
deformation of its classical counterpart (assuming $S$ to be of Hecke type).
At least a necessary condition for flatness of a deformation, i.e., the
existence of a P.b. as a quasiclassical term of the deformation, is
satisfied. The mentioned bracket is also quadratic and admits a
linearization. These two brackets are also compatible. We get a quantization
of the P.p. generated by them by applying to the above algebra a shift
operator (for details the reader is referred to \cite{I}, \cite{IP},
\cite{G2}).}
\end{remark}
\section{Rank 1 quantum orbits}
In what follows we consider the P.p. and their quantum counterparts connected
to symmetric orbits. These P.p. were described in the Introduction: they are
generated by the KKS bracket and by the so-called R-matrix bracket.

We want to describe the corresponding quantum a.p. explicitly by means of a
system of equations. In the present section we recall the method suggested in
\cite{DG1} to look for such systems. Unfortunately, it is valid only for
rank 1 symmetric spaces, namely, those of $\sss$ type.

Let us fix a Lie algebra $g=sl(n)$ and an element $x=\om\in{h}^*$,
where $h\subset g=sl(n)$ is the Cartan subalgebra
(a triangular decomposition of $g$ is assumed to be fixed). We consider
$\om$ as an element of ${g}^*$ extending it by zero to the nilpotent
subalgebras. Let $\Co$ be its orbit in ${g}^*$.

It is well known that this orbit is symmetric iff $\om(h_i)=0$ for all $i\not
=i_0$ where $\{h_i=e_{i,i}-e_{i+1,i+1},\,i=0,....,n-1\}$ is the standard
basis in $h$. The rank of this symmetric space is $rk(\Co)=\min(i_0,n-
i_0)$. We assume here that $rk(\Co)=1$. So we have $i_0=1$ or $i_0=n-1$. To
the sake of concreteness we set $i_0=1$.

Thus, we have a family $\Co$ of such orbits parametrized by the value $\om
(h_1)\in k$ ($k={\bf C}$ or $k={\bf R}$ corresponding to the case under
question). Let us represent the algebra ${\rm Fun}(\Co)$ as a quotient
algebra of ${\rm Sym}(g)={\rm Fun}({g}^*)$. It is not difficult to show
that this algebra can be described by a system of equations quadratic in
generators of the space $g$.

Let us describe this system explicitly. Consider the space $V=g$ as an
object of the category $g$-Mod. Let $\vv=\oplus \vb $ be a decomposition
of $\vv$ into a direct sum of isotypic $g$-modules, where $\beta$ denotes
the highest weight (h.w.) of the corresponding module. Let $\alpha$ be h.w.
of $g$ as a $g$-module.

Then in the above decomposition there is a trivial module $V_0$ ($\beta=0$),
a module $V_{2\alpha}$ of the h.w. $\beta=2\alpha$, and two irreducible
modules isomorphic to $g$ itself. One of them belongs to the space $I_+
\subset \vv$ of symmetric tensors and the other one to the space $I_-\subset
\vv$ of skew-symmetric tensors. We will use the notation $V_+\,(V_-)$ for the
former (latter) of them.

Note that a similar decomposition is valid for other simple Lie algebras,
but for them there exists only one module isomorphic to $g$. It belongs
to the skew-symmetric part and we keep the notation $V_-$ for it.

Let us denote by $C$ the generator (called {\em Casimir}) of the module $V_0$.
Then the orbit $\Co$ is given by the following system of equations
\begin{equation}
I_-=0; \vb=0\,\,\forall\,\vb\subset I_+,\,\,
\vb\not\in\{V_0,\,V_+,\,V_{2\alpha}\};\;C=c_0,\, V_+=c_1 g
\end{equation}
with some factors $c_0$
and $c_1$ (the latter relation is a symbolic way to write the system which
arises if we choose some bases in $g$ and $V_+$ and equate the highest
weight element of $V_+$ to that one of $g$ times a factor $c_1$ and
consider all the relations that follow).

It is not difficult to find the values of the factors $c_i,\,i=0,1$. It
suffices to substitute ``the point'' $\om$ to this system. (Let us note that
in \cite{DG1} this system was given in a inconsistent form, where the
condition $c_1=0$ was assumed.)

A similar procedure can be realized in the category $\uq$-Mod. More
precisely, we equip the space $V$ with the structure of a $\uq$-module
deforming the initial $g$-module structure on $V$(there exists a regular
way to convert any finite-dimensional $g$-module into a $\uq$ one, cf.,
i.e., \cite{CP}). Using the comultiplication in the quantum algebra $\uq$, we
equip the space $\vv$ with an $\uq$-module structure and decompose it as
above into a direct sum of irreducible $\uq$-modules. We will denote these
components by $\vqb$. By $C_q$ we denote a generator of the module $V^q_0$
(we call $C_q$ a {\em braided Casimir}).

Then we can impose a similar system of equations by replacing all $\vb$
participating in the system (4) by their quantum (or q-) analogues $\vqb$ and
$C$ by $C_q$. The only problem is: what are the proper values of the factors
$c_i(q), \, i=0,1$, which now depend on $q$.

In \cite{DG1} the following way to find out a proper system of equations was
suggested. Let us consider the data $(V,\,I^q, \,\nu_0,\,\nu_1)$, where $I^q=
\vv\setminus V^q_{2\alpha}$ and $\nu_i,\, i=0,1$ are two $\uq$-morphisms
defined as follows $\nu_0:V^q_0\to k\;(\nu_0(C_q)=c_0)$ and $\nu_1:V^q_+\to
c_1 g$ (all other components of $\vv\setminus V^q_{2\alpha}$ are sent by
$\nu_i$ to zero).

Let us also consider the graded quadratic algebra $T(V)/\{I^q\}$ and its
filtered analogue $T(V)/\{\iqc\}$, where $\{\iqc\}$ the ideal generated by
the elements
$$
I^q_-,\, \vb\subset I_+^q,\,\,
\vb\not\in\{V_0,\,V_+,\,V_{2\alpha}\},\;C-c_0(q),\, V_+^q-c_1(q)g
$$
(see the end of this Section for the definition of the spaces $I^q_{\pm}$).

Let us note that the algebra $T(V)/\{I^q\}$ is the quantum analogue of the
function space on the cone and it is a flat deformation of its classical
counterpart. The flatness can be shown, e.g., in the way suggested in
\cite{DG2}, where an intertwining operator of the initial commutative product
and the deformed one is given. Moreover, the algebra $T(V)/\{I^q\}$ is Koszul
for a generic $q$ since it is so for the case $q=1$ by virtue of \cite{B}.

Let us assume now that the above data satisfy the following system
\begin{eqnarray*}
&&{\rm Im}(\nu_1\ot{\rm id}-{\rm id}\ot\nu_1)(I\ot V\bigcap V\ot I)\subset
I,\\
&&(\nu_1(\nu_1\ot {\rm id}- {\rm id} \ot \nu_1)+
\nu_0 \ot {\rm id} - {\rm id} \ot \nu_0)(I \ot V\bigcap V\ot I)=0,\\
&&\nu_0 (\nu_1\ot {\rm id}- {\rm id} \ot \nu_1)(I \ot V\bigcap V\ot I)=0.
\end{eqnarray*}
Then by virtue of the PBW theorem in the form of \cite{BG} we can conclude
that the algebra $T(V)/\{\iqc\}$ is a flat deformation of the algebra $T(V)/
\{I^q\}$.

Let us note that the above conditions represent a more general form of the
Jacobi relation connected to deformation theory. Thus, if the above form of
the Jacobi identity is fulfilled, the algebra $T(V)/\{\iqc\}$ is a flat
deformation of the orbit $\Co$ (more precisely, of the function algebra on
it).

So, the proper quantities $c_i(q),\,i=0,1$, if they exist, can be found from
the above equations. However, a priori it is not clear why such quantities
exist. Let us assume here that they exist and denote by $\aq$ the
corresponding algebra $T(V)/\{\iqc\}$. We are interested in its further
deformation.

To do this, we begin by discussing the following question: what is the
deformational quantization of the KKS bracket? We want to represent the
latter quantum object also as a quotient of the enveloping algebra $U(sl(n))
_h$ The index $h$ here means that we have introduced a factor $h$ in
the bracket in the definition of the enveloping algebra, i.e., $$U(g)_h=
T(g)/\{xy-yx-h[x,y]\}.$$

There are some factors $c_i(h),\,i=0,1$, now depending on $h$ such that the
quotient algebra $\ah=U(sl(n))_h/\{J\}$, where the ideal $\{J\}$ is
generated by the family of elements from (4) lying in $I_+$ but with new
$c_i(h)$, is a flat deformation of the initial algebra corresponding to the
case $h=0$. These factors can be also found by means of the above form of
Jacobi identity.

This approach can also be applied in order to get the quantum algebras
corresponding to the whole P.p. under consideration, since these algebras are
quadratic as well. The only problem is: what is the proper quantum analogue
of the algebra $U(sl(n))_h$? Or, in other words, what is a consistent way
to introduce a quantum analogue of the Lie bracket?

We will introduce a q-generalization of the ordinary Lie bracket by means of
the following
\begin{definition}{\em
Let $g$ be a simple Lie algebra
equipped with a $\uq$-module structure and $\vqb$ be the irreducible modules
in the category $\uq$-Mod entering its tensor square. We call {\em q-Lie
bracket} the operator $[\,\,,\,\,]_q:\vv\to V$ defined as follows
$$[\,\,,\,\,]_q|_{\vb^q}=0\mbox{ for all }\vb^q\not=V_-^q$$
and $ [\,\,,\,\,]_q:V_-^q\to V$ is a $\uq$-isomorphism.}
\end{definition}

Let us observe that in this way the q-Lie bracket is defined up to a factor.
We fix this factor and introduce the q-analogue of the algebra $U(g)_h$
as follows
$$U(g)_{h,q}=T(V)/\{{\rm Im}({\rm id}-h[\,\,,\,\,]_q)I_-^q\}.$$
The space $V=g$ equipped with this bracket will be denoted by $\ogg$ and
called {\em braided Lie algebra}\footnote{Other definitions of q-counterparts
of the Lie algebras have been suggested recently in \cite{DH} and \cite{LS}.
It is very plausible that they are equivalent to ours.}. We will also use the
notation $\uog$ for the algebras $U(g)_{1,q}$. Let us note that this is
another, as compared with the quantum group $\uq$, q-analogue of the
enveloping algebra $U(g)$, but the deformation $U(g)\to U(g)_{h,
q}$ is not flat except for the $sl(2)$ case. However, some quotient algebras
of this algebra are flat deformations of their classical counterparts.

Let us return to the case $g=sl(n)$ and introduce the first type a.p.
$\ahq$ as the quotient algebra of $U(sl(n))_{h,q}$ by the ideal
$$\{C_q-c_0,\,V_+^q-c_1V,\,I^q_+\setminus(V^q_{2\alpha}\oplus V_0^q\oplus
V_+^q)\}$$
with some factors $c_i(h,q),\,i=0,1$.

In order to look for the consistent factors $c_i(h,q)$, we must only modify
the above morphism $\nu_1$ on the ``q-skew-symmetric'' subspace $I^q_-$ by
setting $\nu_1: V_-^q\to h g$ (all other components are still sent by $\nu
_1$ to zero) and verify the above form of Jacobi identity. This provides us
with the factors $c_i(h,q)$ ensuring the flatness of the deformation of the
function algebra on the orbit $\Co$ under consideration.

Let us emphasize that the existence of the proper factors $c_i(h,q)$, as well
as that of the above factors $c_i(q)$, can be deduced from the paper
\cite{DS1}, where a formal quantization of the P.p. under question was
considered. In the next Section we discuss another way to look for appropriate
factors $c_i(h,q)$.

It remains only to note that the quasiclassical term of two-parameter
deformation ${\rm Fun}(\Co) \to \ahq$ is just the above P.p. on the orbit
$\Co$ (cf. \cite{DG1}).

It is worth to note that the product in the resulting quantum algebras $\ahq$
is $\uqsn$-invariant. This means that the following property
$$X\mu(a\ot b)=\mu(X_1(a)\ot X_2(b)),\; a,b\in\ahq,\; X\in\uqsn$$
is satisfied. (Hereafter we use Sweedler's notation $X_1\ot X_2$ for $\De
(X)$.) This property follows immediately from the construction of the
algebras. This algebra is the inductive limit of finite-dimensional
$\uqsl$-modules (moreover, it is multiplicity free, i.e., the multiplicity of
each irreducible $\uq$-module is $\leq 1$). Thus, this algebra belongs to the
category $\uqsl$-Mod.

Let us note that a particular case of the a.p. $\ahq$ is the algebra $\aq$
arising as the result of quantization of the only R-matrix bracket. This
algebra is commutative in the category $\uqsl$-Mod in the following sense.
This category is balanced (see \cite{CP} for the definition). Moreover, for
any two finite-dimensional objects $U$ and $V$ of this category there exists
an involutive operator $\tS:U\ot V\to V\ot U$ such that it is a
$\uq$-morphism, it commutes with $S$ and it is a deformation of the flip.

Then $I^q_{\pm}={\rm Im}({\rm id}\pm\tS)$ with $\tS:\vv\to\vv$. Using the fact
that the algebra $\ahq$ can be decomposed into a direct sum of finite-dimensional
objects of the category $\uqsl$-Mod we can extend the operator $\tS$ onto
$\ahq^{\ot 2}$. Then the above mentioned commutativity of the algebra $\aq$
means that the multiplication operator $\mu$ satisfies the following relation
$\mu=\mu \tS$. A proof of this fact can be obtained from \cite{DS1}.
Strictly speaking, just the algebra $\aq$ is the quantum analogue of the
orbit under consideration.

Let us note that the above method can be also applied to find the equations
describing quantum algebras arising from similar P.p. on certain nilpotent
orbits in ${g}^*$ (see the Introduction). Thus, such a.p. related to the
highest weight element orbits in ${g}^*$ can be defined by means of the
above ideal $\iqc$, but with $c_0(h,q)=0$ and with a suitable $c_1(h,q)$.

\section{Modules for first type quantum algebras}
In the previous section we discussed the first step of quantization
procedure. The resulting object of this step is an a.p. represented as a
quotient algebras over some suitable ideal. Now we want to consider the
second step of quantization, consisting in an attempt to represent the above
algebras in certain linear spaces. At this step the difference between two
types of quantum algebras (a.p.) under question becomes clearer.

Moreover, this step provides us with another way to look for consistent
factors $c_i(h,q)$. Briefly speaking, this method reduces to the computation
of the factor $c_i(h,q)$ on the image of the first type algebras into the
space ${\rm End}(V)$.

To do this, we recall a natural way to equip the space of endomorphisms of a
$\uq$-module with a $\uq$-module structure. Let $U^q$ be a
(finite-dimensional) $\uq$-module and $\rho_q:\uq\to {\rm End}(U^q)$ be the
corresponding representation. Let us introduce the representation $\ren:\uq
\to{\rm End}({\rm End}(U^q))$ by putting
$$
\ren(a)M=\rho(a_1)\circ M\circ\rho(\gamma(a_2)),\,a\in\uq, \,M\in {\rm End}
(U^q).
$$
We denote the matrix product by $\circ$, while $\gamma$ is the antipode in
$\uq$.

We deal with a coordinate representation of module elements. We consider the
endomorphisms as matrices and their action is the left multiplication by
these matrices.

Let us note that this way of equipping ${\rm End}(U^q)$ with a $\uq$-module
structure is compatible with the matrix product in it in the following
sense:
$$\ren(a)(M_1\circ M_2)=\ren(a_1)M_1\circ\ren(a_2)M_2.$$
This means that $\circ: {\rm End}(U^q)^{\ot 2}\to {\rm End}(U^q)$ is a
$\uq$-morphism.

Now we will introduce a useful notion for constructing a representation
theory of the algebras under consideration.

\begin{definition}{\em Let $\ogg$ be a braided Lie algebra.
We say that a $\uq$-module $U^q$ is a {\em braided module} or, more
precisely, a {\em braided $\,\ogg$-module} if it can be equipped with a
structure of a $\uog$-module and the representation $\rho:U(\ogg)\to {\rm
End}(U^q)$ is a $\uq$-morphism. We also say that the classical counterpart
$U=U^1$ of the $\uq$-module $U^q$ allows braiding. }
\end{definition}

A natural way to construct braided modules is given by the following

\begin{proposition} Let $U^q$ be a $\uq$-module. If the decomposition ${\rm
End}(U^q)=\oplus V^q_{\gamma}$ of the $\uq$-module ${\rm End}(U^q)$ into the
direct sum of irreducible $\uq$-modules is such that

{\rm 1.} it does not contain modules isomorphic to $\vbq\subset I_-^q$
apart from those isomorphic to $V^q$, where by $V^q$ we denote $g=V$
equipped with the $\uq$-module structure,

{\rm 2.} the multiplicity of the module $V^q$ is $1$,

\noindent then $U^q$ can be equipped with a braided module structure
{\rm (}briefly, braided structure{\rm )}.
\end{proposition}

{\it Proof}. By the assumption, there exists a unique $\uq$-submodule in
${\rm End}(U^q),$ isomorphic to $V^q$. Consider a $\uq$-morphism defined up
to a factor
\begin{equation}
\rho: V^q\to {\rm End}(U^q).
\end{equation}
This map is an almost representation of the braided Lie algebra $\ogg$ in the
sense of the following
\begin{definition}{\em We say that a map (5) is an {\em
almost representation} of the braided Lie algebra $\ogg$ if it is a
$\uq$-morphism and the following properties are satisfied

1. $\circ\rd(\vbq)=0$ for all $\vbq\subset I_-^q$ apart from that
$V^q_-$,

2. $\circ\rd V_-^q=\nu\rho[\,\,,\,\,]_qV^q_-$ with some $\nu\not=0$.}
\end{definition}
In a more explicit form these conditions can be reformulated as follows. If
the elements $\{b^{i,j}_{k,\beta}u_iu_j, 1\leq k\leq\dim\, \vbq\}$ form a
basis of the space $\vbq\subset I_-^q$ and similarly the elements $\{b^{i,
j}_{k,-}u_iu_j, 1\leq k\leq\dim\, V^q_-\}$ form a basis of the space $V_-^q$
then
$$
b^{i,j}_{k,\beta}\rho(u_i)\rho(u_j)=0\,\,{\rm if}\, \,\vbq\not=V_-^q\, \,{\rm
and}\,\, b^{i,j}_{k,-}\rho(u_i)\rho(u_j)=\nu b^{i,j}_{k,-}\rho([u_i,u_j]_q).
$$

Let us complete the proof. The image of the composed map $\circ\rd(V^q_-)$ is
isomorphic to the $\uq$-module $V^q$ since $\rho$ and $\circ$ are
$\uq$-morphisms. Such a module in the decomposition of ${\rm End}(U^q)$ is
unique, therefore the image of the space $[\,\,,\,\,]_q V^q_-=V^q$ with
respect to the morphism $\rho$ coincides with the previous one (it suffices
to show that the images of the highest weight element of the module $V_-^q$
with respect to both operators coincide up to a factor). This gives the
second property of Definition 2 (the property that $\nu\not=0$ for a generic
$q$ follows from the fact that this is so for $q=1$).

The first property of the Definition follows from the fact that ${\rm End}
(U^q)$ does not contain any modules isomorphic to $\vbq\subset I^q_-,\, \vbq
\not=V_-^q$. Finally, changing the scale, i.e., considering the map $\rho_
{\nu}=\nu^{-1}\rho$ instead of $\rho$, we get a representation of the algebra
$\uog$. This completes the proof. \ $\Box$

A natural question arises: how many $\uq$-modules can be converted into
braided ones or, in other words, how many $g$-modules allow braiding? The
answer to this question for the $sl(n)$-case is given by the following
proposition, which can be proved by straightforward computations using Young
diagram techniques.

\begin{proposition} Let $\omega$ be a fundamental weight of the Lie algebra
$sl(n)$. Then the $sl(n)$-modules $\vko$ {\em (}for any nonnegative
integer $k${\em )} allow  braiding. In other words, their q-analogues
$\vko^q$ are braided modules.
\end{proposition}

For other simple Lie algebras, $g$ it seems very plausible that a similar
statement is valid for fundamental weights $\omega$ such that their orbits in
${g}^*$ are symmetric (in the $sl(n)$-case all fundamental weights
satisfy this condition). Note that all orbits of such type have been
classified by E. Cartan (cf. \cite{KRR}).

Let us now discuss how the braided modules can be used to find the above
factors $c_i(h,q)$.

Once more we set $g=sl(n)$ and we take as fundamental the weight $\om$
such that one $\om(h_1)=1,\, \om(h_i)=0,\, i>1$. Then we consider the $U(sl
(n))$-module $\vkoq$ (which is in fact $\vko$ equipped with a
representation $\rho_q :\uq \to {\rm End}(\vko^q)$). Let us realize a
braiding of the modules, i.e., construct a $U_q(g)$-morphism
$$\rho:U(g)_{h,q}\to {\rm End}(\vkoq)$$
in the way described above for the algebra $U(\ogg)=U(g)_{1,q}$.

\begin{proposition} The representation $\rho$ is factorized to a
representation of the algebras $\ahq$ with certain $c_i(h,q)$. {\em (}This
means that $\rho(\vbq)=0$ if $\vbq\subset I^q_+\setminus(V^q_{2\alpha}
\oplus V_0^q \oplus V_+^q)\}$, $\rho(C_q)=c_0(h,q) {\rm id}$, and $\rho(V_+
^q)=c_1(h,q)\rho(V^q)$.{\em )}
\end{proposition}

{\it Proof.} It suffices to show that the modules belonging to $I^q_+
\setminus(V^q_{2\alpha}\oplus V_0^q\oplus V_+^q)\}$ are not represented in
${\rm End}(\vko^q)$ for any $k$ and that multiplicity of the modules $V^q_0$
and $V_+^q$ in ${\rm End}(\vko^q)$ is 1. This can be done by straightforward
calculations by means of Young diagram techniques.

Finally, the factors $c_i(h,q)$ are defined by these relations. Now if we
want to find the relations defining the ``quantum orbits'' (i.e. algebras
corresponding to the case $h=0$) we pass to the limits $k\to\infty$ and $h\to
0$ in such a way that $c_0(h,q)$ has a limit (denoted by $c_0(q)$). By virtue
of \cite{DS1} $c_1(h,q)$ also has a limit (denoted by $c_1(q)$). These two
constants are just the factors that we are looking for. \ $\Box$

It is interesting to emphasize that to find the system of equations
describing a commutative algebra in the category $\uq$-Mod in the framework
of this approach, we looking first for the system corresponding to the
noncommutative algebra $\ahq \,(h\not=0)$.

\section{On ``twisted quantum mechanics''}

In this section we discuss the so-called twisted quantum mechanics. Roughly
speaking, it is a quantum mechanics in twisted categories. Our aim is to
consider the above quantum algebras from this point of view.

It is well known (De Wilde--Lecomte--Fedosov) that any symplectic Poisson
structure can be quantized and the result of quantization can be treated as
an operator algebra ${\rm End}(V)$ in a complex Hilbert space $V$. The space
${\rm End}(V)$ is equipped in this case with a trace and a conjugation
(involution) possessing the usual properties
$${\rm tr}(A\circ B)={\rm tr}(B\circ A),\;{\rm tr}\,A^*={\rm tr}\, A,\;
(A\circ B)^*=B^*\circ A^*.$$

The quantum observables are identified with Hermitian (self-adjoint)
operators in this space and they form a linear space over ${\bf R}$ closed
with respect to the bracket $i[\,\,,\,\,]$.

If we consider a super-version of quantum mechanics, this bracket must be
replaced by its super-analogue. The notions of ordinary trace and
conjugation operators must be also replaced by their super-counterparts; this
leads, in particular, to a modification of the notion of self-adjoint
operators. Meanwhile, the observables that play the role of Hamiltonians
must be even (since only for an even operator $H$ the equation $dA/dt=i[H,A]$
is consistent) and such an operator is self-adjoint in the usual sense iff it
is super-self-adjoint. Moreover, the super-bracket with such an operator
becomes the ordinary Lie bracket.

What is a proper analogue of this scheme in a twisted category? If a
category under question is symmetric, i.e., the Yang--Baxter twist $S$ is
involutive $(S^2={\rm id})$ in it, the corresponding generalization was
suggested in \cite{GRZ}. In this case ``S-analogues'' of Lie bracket, trace
and conjugation operators introduced there satisfy the following version of
the above relations
\begin{equation}
{\rm tr}(A\circ B)={\rm tr}\circ S(A\ot B),\;{\rm tr}\,A^*={\rm tr}\, A,\;
(A\circ B)^*=\circ(*\ot *)S(A\ot B).
\end{equation}

Let us recall a construction of a conjugation operator in this case (assuming
$V$ to be finite-dimensional). Let us suppose that we can identify $V$ and
$V^*$ (this means that there exists a pairing $\langle\,\,,\,\,\rangle:\vv\to
k$, which is a morphism of the category). Then ${\rm End}(V)$ can be
identified with $\vv$ and the conjugation $*:{\rm End}(V)\to{\rm End}(V)$ is
just the image of the operator $S:\vv\to \vv$ under this identification.
Therefore the conjugation satisfies the following relation $(*\ot {\rm id})S=
S({\rm id}\ot *)$.

More precisely, we consider the space $g=V$ over the field $k={\bf R}$
and assume all matrix elements of $S$ to be real. This implies that if we
introduce an S-Lie bracket in ${\rm End}(V)$ by
$$
[A,B]=A\circ B-\circ S(A\ot B)
$$
its structural constants are real as well.

Let us extend this bracket to the space $V_{{\bf C}}=V\ot{\bf C}$ by
linearity. Let $*: V_{{\bf C}}\to V_{{\bf C}}$ be a conjugation, i.e., an
involutive $(*^2={\rm id})$ operator such that $(\lambda z)^*=\overline{
\lambda}z^*,\,\lambda\in {\bf C},\, z\in V_{{\bf C}}$. 
Assume also that the relation $(*\ot {\rm id})S=S({\rm id}\ot *)$ and the third
relation (6) are satisfied. Then it is easy to
see that this conjugation is compatible with the S-Lie bracket in the
following sense
\begin{equation}
*[A,B]=-[A^*,B^*].
\end{equation}

It is not so evident what is the proper definition of a conjugation operator,
say, in the algebra ${\rm End}(V)$, where $V$ is an object of a twisted but
nonsymmetric category. Now possessing a q-analogue of the Lie bracket, we can
try to define the compatibility of such an operator with the twisted
structure by means of relation (7). (Let us note that the operator $S$ does
not enter explicitly in it and formula (7) is in some sense universal.)

Let us consider a space $V_{{\bf C}}$ equipped with a bracket
$[\,\,,\,\,]:V_{{\bf C}}^{\ot 2}\to V_{{\bf C}}$.

\begin{definition}{\em We say that a conjugation $*$ is {\em compatible}
with this bracket  if the relation (7) is satisfied.}
\end{definition}
\begin{proposition} The odd elements with respect to this involution
{\em (}i.e., elements such that $z^*=-z${\em )} form a subalgebra, i.e.,
the element $[a,\,b]$ is odd if $a$ and $b$ are.
\end{proposition}

{\it Proof}. It is obvious. \ $\Box$

Therefore the space of ``$*$-even'' operators is closed with respect to the
bracket $i[\,\,,\,\,]$.

In \cite{DGR}, we have classified all such conjugations for $\uqs$-case. Let
us reproduce the final result here, but first represent the multiplication
table for q-Lie bracket in some base $\{u,v,w,\}$:
$$[u,u]=0,\ [u,v]=-q^2Mu,\ [u,w]=(\qq)^{-1}Mv,$$
$$[v,u]=Mu,\ [v,v]=(1-q^2)Mv,\ [v,w]=-q^2Mw,$$
$$[w,u]=-(\qq)^{-1}Mv,\ [w,v]=Mw,\ [w,w]=0,$$
Here $M$ is an arbitrary real factor (cf. \cite{G2} for details).

\begin{proposition} For a real $q\not=1$ there exist only two conjugations in
the space $V_{{\bf C}}$ compatible with q-Lie bracket, namely, that $a^*=
-\overline a$ for any $a\in V_{{\bf C}}$ and the following one $u^*=u,\,
v^*=-v,\,w^*=w$. \end{proposition}

Although these conjugations are rather trivial they are, together with
quantum traces (their construction in ${\rm End}(V)$, where $V$ is an object
of a rigid category, is well known) the ingredients of twisted quantum
mechanics in the sense of the following informal definition.

\begin{definition}{\em
We say that an associative algebra is a {\em subject of twisted quantum
mechanics} if it belongs to a twisted category, it is represented in the
space ${\rm End}(V)$ equipped with a twisted Lie bracket, a trace and a
conjugation as above and the representation map is a morphism in this
category.}\end{definition}

Unfortunately, we cannot give a final axiom system of twisted quantum
mechanics for a noninvolutive $S$ (though it seems very plausible that the
above relations between the operators under consideration are still valid in
the category $\uq$-Mod if we replace $S$ by $\tS$). However, we want to point
out the principal difference between the twisted version of quantum mechanics
and its classical version: the twisted (quantum) trace must occur in
calculations of partition functions.

In the above sense the second type quantum algebras represented in some linear
spaces in spirit of the paper \cite{VS} are not subjects of
twisted quantum mechanics. Though the algebra $\aq={\rm Sym}(W)$ itself
belongs to the twisted category generated by the space $W$, its modules 
constructed  in \cite{VS} for $su(2)$ case do not belong to this category.

As for the second type algebras $\ahq$ with $h\not=0$, they differ
nonessentially from the algebras $\aq$.

This is not the case for the first type algebra $\ahq$. The representation
theory of the algebras $\ahq$ for generic $h$ and that for the case $h=0$ are
completely different. From this point of view $h=0$ is a singular point and we
disregard it.

Let us note that our first type quantum algebras are subjects of twisted
quantum mechanics arising from the quantization of Poisson brackets. However,
it is possible to consider similar objects connected to nondeformational
solutions of the QYBE (cf. \cite{G1}, \cite{GRZ}).

Concluding the paper, we would like to formulate two problems: that of
calculating the above factors $c_i(h,q)$ (or more generally, giving an exact
description of the two parameter deformation of all symmetric orbits in
${g}^*$ for any simple Lie algebra $g$) and the problem of generalizing
our approach to infinite-dimensional Lie algebras.


\begin{thebibliography}{DGM}


\bibitem[B]{B} R.B\"orgvad 
{\em Some homogeneous coordinate rings that are Koszul algebras}, alg-geom/951011.


\bibitem[BG]{BG} A.Braverman, D.Gaitsgory
{\em Poincar\'e-Birkhoff-Witt theorem for quadratic algebras
of Koszul type}, hep-th/9411113.

\bibitem[CP]{CP} V.Chari ,
A.Pressley   {\em A guide to Quantum Groups},  Cambrige University press, 1994

\bibitem[DG1]{DG1} J.Donin, D.Gurevich {\em Quantum orbits of R-matrix
type}, Lett. Math. Phys. 35 (1995), pp. 263--276.

\bibitem[DG2]{DG2}
J.Donin, D.Gurevich {\em Some Poisson structures associated to Drinfeld-Jimbo R-matrices
and their quantization}, Israel Math. Journal 92 (1995), pp.23-32                                                                       
 
\bibitem[DGM]{DGM}
J.Donin , D.Gurevich, S. Majid  {\em R-matrix brackets and their quantization}
Ann. Inst. Henri Poincar\'e  58 (1993), pp. 235-246


\bibitem[DGR]{DGR} J.Donin, D.Gurevich, V.Rubtsov  {\em Quantum hyperboloid and
braided modules}, Proceedings of French-Belgium meeting, 1995, Reims,  to
appear

\bibitem[DS1]{DS1} 
J.Donin , S.Shnider {\em Quantum symmetric spaces} 
 Journal of pure and applied algebra 100 (1995), pp. 103-115
 
\bibitem[DS2]{DS2} 
J. Donin, S.Shnider  {\em Quasi-associativity and flatness criteria for 
quadratic algebra deformation},  Israel Math. J., to appear
 
\bibitem[G1]{G1} D.Gurevich {\em Algebraic aspects of the quantum 
Yang-Baxter equation}, Leningrad Math.J. 2 (1991),  
pp.801--828.

\bibitem[G2]{G2}
Gurevich D. {\em Braided modules and reflection equations}, Publications of Banach Center, 
Warsaw, to appear

\bibitem[GP]{GP} D.Gurevich, D.Panyushev {\em On Poisson pairs associated to
 modified R-matrices}, Duke Math. J. 73 (1994), pp.249--255.

\bibitem[GR]{GR} D.Gurevich, V.Rubtsov, {\em Quantization of Poisson pencils and
generalized Lie algebras}, Teor. i Mat. Phys. 103 (1995), pp. 476--488.

\bibitem[GRZ]
{GRZ} D.Gurevich , V.Rubtsov, N.Zobin {\em Quantization of 
Poisson pairs: R-matrix approach},  Journ.Geom.and Phys. 9 (1992), pp.25--44

\bibitem[I]{I} A.Isaev {\em Interrelation between Quantum Groups and
Reflection Equation (Braided) Algebras}, Lett. Math. Phys. 34 (1995), pp.
333--341.

\bibitem[IP]{IP} A.Isaev, P.Pyatov {\em Covariant differential complex on
quantum linear groups}, J. Phys. A: Math. Gen. 28 (1995), pp. 2227--2246.

\bibitem[KRR]{KRR}
S.Khoroshkin, A.Radul, V.Rubtsov{\em Families of Poisson structures on
Hermitean Symmetric Spases}, Comm. Math. Phys. 153 (1993), pp.299-315

\bibitem[RTF]{RTF} 
N.Reshetikhin, L.Takhtadzhyan, L. Faddeev {\em Quantization
 of Lie groups and Lie algebras}, Leningrad Math.J. 1 (1990),
pp.S193--226

\bibitem[S-T]{S-T}
M. Semenov-Tian-Shansky {\em What is classical r-matrix?},
  Funct. Anal. Appl. 17 (1883),  pp.S259--272

\bibitem[S]{S} E.Sklyanin 
 {\em Some algebraic structures connected with the Yang-Baxter equation}
 Funct. Anal. Appl. 16 (1982), pp.S263-270

\bibitem[VS]{VS}
L.Vaksman, Y.Soibelman  {\em Algebra of functions on the quantum qroup $su(2)$},
Funct. Anal. Appl. 22 (1988),  pp.S170-181

\end{thebibliography}
 \end{document}